\documentstyle[12pt]{article}
\input amssym.def
\begin{document}
\def \Z{\Bbb Z}
\def \C{\Bbb C}
\def \R{\Bbb R}
\def \Q{\Bbb Q}
\def \N{\Bbb N}
\def \wt{{\rm wt}}
\def \tr{{\rm tr}}
\def \span{{\rm span}}
\def \Res{{\rm Res}}
\def \Res{{\rm QRes}}
\def \End{{\rm End}}
\def \E{{\rm End}}
\def \Ind {{\rm Ind}}
\def \Irr {{\rm Irr}}
\def \Aut{{\rm Aut}}
\def \Hom{{\rm Hom}}
\def \mod{{\rm mod}}
\def \ann{{\rm Ann}}
\def \<{\langle} 
\def \>{\rangle} 
\def \t{\tau }
\def \a{\alpha }
\def \e{\epsilon }
\def \l{\lambda }
\def \L{\Lambda }
\def \g{\gamma}
\def \b{\beta }
\def \om{\omega }
\def \o{\omega }
\def \c{\chi}
\def \ch{\chi}
\def \cg{\chi_g}
\def \ag{\alpha_g}
\def \ah{\alpha_h}
\def \ph{\psi_h}

\def \bconj{\begin{conj}\label}
\def \econj{\end{conj}}
\def \be{\begin{equation}\label}
\def \ee{\end{equation}}
\def \bex{\begin{exa}\label}
\def \eex{\end{exa}}
\def \bl{\begin{lem}\label}
\def \el{\end{lem}}
\def \bt{\begin{thm}\label}
\def \et{\end{thm}}
\def \bp{\begin{prop}\label}
\def \ep{\end{prop}}
\def \br{\begin{rem}\label}
\def \er{\end{rem}}
\def \bc{\begin{coro}\label}
\def \ec{\end{coro}}
\def \bd{\begin{de}\label}
\def \ed{\end{de}}
\def \pf{{\bf Proof. }}
\def \voa{{vertex operator algebra}}

\newtheorem{thm}{Theorem}[section]
\newtheorem{prop}[thm]{Proposition}
\newtheorem{coro}[thm]{Corollary}
\newtheorem{conj}[thm]{Conjecture}
\newtheorem{exa}[thm]{Example}
\newtheorem{lem}[thm]{Lemma}
\newtheorem{rem}[thm]{Remark}
\newtheorem{de}[thm]{Definition}
\newtheorem{hy}[thm]{Hypothesis}
\makeatletter
\@addtoreset{equation}{section}
\def\theequation{\thesection.\arabic{equation}}
\makeatother
\makeatletter

\newcommand{\rw}{\rightarrow}
\newcommand{\n}{\:^{\circ}_{\circ}\:}

\begin{center}{\Large \bf Some finiteness properties of regular
vertex operator algebras}\\
\vspace{0.5cm}
Haisheng Li\footnote{Partially supported by NSF grant DMS-9616630.}\\ 
Department of Mathematical Sciences\\
Rutgers University-Camden\\
Camden, NJ 08102
\end{center}

\begin{abstract}
We give a natural extension of the notion of the contragredient module
for a vertex operator algebra. By using this extension we prove that 
for regular vertex operator algebras, Zhu's $C_{2}$-finiteness condition
holds, fusion rules (for any three irreducible modules) are finite and 
the vertex operator algebras themselves are finitely generated.
\end{abstract}

\section{Introduction}
In this paper, we study certain finiteness properties of 
vertex operator algebras, motivated by 
a conjecture of Zhu about the relation between the rationality and 
the $C_{2}$-finiteness condition for a vertex operator algebra $V$.

To prove the convergence of a trace function
Zhu [Z] made a technical assumption on the vertex operator algebra, 
called finiteness condition $C$,
so that his beautiful results hold for a rational vertex operator algebra
that satisfies the finiteness condition $C$. 
Zhu's finiteness condition $C$ consists of what we call
$C_{2}$-finiteness condition in the present paper and the condition that
$V$ is a sum of lowest weight modules for the Virasoro algebra.
We say that a vertex operator algebra $V$ satisfies {\em $C_{2}$-finiteness 
condition} if $C_{2}(V)$ is finite-codimensional in $V$ 
where $C_{2}(V)$ is the 
subspace of $V$ linearly spanned by elements $u_{-2}v$ for $u,v\in V$.
(Note that as one of the results in [DLM3], the convergence of trace 
functions was proved under only the $C_{2}$-finiteness condition.)
It was proved in [Z] (see also [DLM3]) that
the familiar rational vertex operator algebras satisfy
the $C_{2}$-finiteness condition and
it was conjectured that the rationality in the sense of [Z] 
(defined in Section 2) implies the $C_{2}$-finiteness condition. 

In [DLM2], with various motivations we proved that for 
the familiar rational vertex operator algebras,
any weak module (defined in Section 2) is a direct sum of 
irreducible (ordinary) modules.
Consequently, any irreducible weak module is an (ordinary) module.
Vertex operator algebras with this property are said to be {\em regular}.
It was conjectured in [DLM2] that
the notions of rationality and regularity are equivalent. The combination of
the two conjectures gives rise to a third conjecture:
Regularity implies the $C_{2}$-finiteness condition.
Since regularity clearly implies rationality, the third conjecture
is also a weak version of Zhu's conjecture.

In this paper, as one of our main results we prove the third conjecture, 
which is that regularity implies the $C_{2}$-finiteness condition.
Moreover, we prove that if $V$ is regular, the fusion rules are finite 
and $V$ is finitely generated as a vertex operator algebra.

%Notice that Zhu's $C_{2}$-finiteness condition does not depend on
%the choice of the Virasoro vector and that roughly speaking
%regularity does not depend on the choice of the Virasoro vector either.
%However, rationality (at least superficially)
%does depend on the choice of the Virasoro vector because
%the $L(0)$-weight grading was used for the notion of ${\N}$-gradable weak 
%module. So philosophically speaking, Zhu's $C_{2}$-finiteness condition is 
%closer to regularity than rationality. 

In [Z], Zhu associated an associative algebra $A(V)$ (a certain quotient 
space of $V$) to any vertex 
operator algebra $V$ and established natural functors between the category of 
${\N}$-gradable weak $V$-modules (defined in Section 2) and 
the category of $A(V)$-modules.
Furthermore, these functors give rise to a one-to-one 
correspondence between the set of equivalence classes of 
irreducible $A(V)$-modules and the set of equivalence classes 
of irreducible ${\N}$-gradable weak $V$-modules. 
Using these functors Zhu was able to prove that
$A(V)$ is (finite-dimensional) semisimple if $V$ is rational.
Although the quotient space $Z_{2}(V)=V/C_{2}$
has a natural commutative associative algebra structure [Z], there are no
appropriate functors between the category of $Z_{2}(V)$-modules and
the category of ${\N}$-gradable weak $V$-modules. Also notice that
if $V$ satisfies Zhu's $C_{2}$-finiteness condition and the condition
$V=\coprod_{n=0}^{\infty}V_{(n)}$ with $V_{(0)}={\C}{\bf 1}$, then 
the nil-radical of
$Z_{2}(V)$ is $1$-codimensional (because $Z_{2}(V)$ is an ${\N}$-graded
algebra). Thus, in principle there are no desired functors between 
the category of $V$-modules and the category of $Z_{2}(V)$-modules.
This indicates that a different technique is needed.

In the following we give an account of the main stream of this paper.
Let $W=\coprod_{h\in {\C}}W_{(h)}$ be a $V$-module. Then we 
generalize the definition of $C_{2}(V)$ to define $C_{2}(W)$
in the obvious way.
View $(W/C_{2}(W))^{*}$ as a natural subspace of $W^{*}$.
Since $C_{2}(W)$ is a graded subspace, $C_{2}(W)$ is finite-codimensional 
if and only if
$(W/C_{2}(W))^{*}\subseteq W'$ $(=\coprod_{h}W_{(h)}^{*}$,
the restricted dual of $W$). 
To use the rationality or the regularity, first of all one needs to
relate $(W/C_{2}(W))^{*}$ to a certain weak $V$-module.

A fundamental result proved in [FHL] is that $(W',Y')$
carries the structure of a $V$-module where
$$\< Y'(v,x)f,w\>=\<f,Y(e^{xL(1)}(-x^{-2})^{L(0)}v,x^{-1})w\>$$
for $v\in V, f\in W', w\in W$.
One notices that
the action of $Y'(v,x)$ on $W'$ obviously extends to $W^{*}$, say $Y^{*}(v,x)$.
But $W^{*}$ fails to be even a weak $V$-module.
If there is a weak $V$-module $M$ in $W^{*}$ containing
$(W/C_{2}(W))^{*}$, on which $L(0)$ acts semisimply, then it is
easy to prove that $M\subseteq W'$, so that
$(W/C_{2}(W))^{*}\subseteq W'$. This would solve our problem.
Notice that $L(0)$ acts semisimply on any weak (${\N}$-gradable weak)
module for a regular (rational) vertex operator algebra.
Now, the question is whether
there exists a weak or an ${\N}$-gradable weak $V$-module in $W^{*}$ containing
$(W/C_{2}(W))^{*}$ as a subspace.
If $S$ and $T$ are subspaces of $W^{*}$ such that $(S,Y^{*})$ and $(T,Y^{*})$ 
are 
weak $V$-modules, then $(S+T,Y^{*})$ is clearly a weak $V$-module again.
Thus there is a unique maximal weak $V$-module in $W^{*}$ with the vertex 
operator map $Y^{*}$.
Then we are naturally led to this weak module to see whether it contains
$(W/C_{2}(W))^{*}$.

Notice that in the notion of module or weak module $W$, the first axiom is 
the truncation condition: $Y(v,x)w\in W((x))$ for $v\in V,w\in W$, and that
without this condition the Jacobi identity cannot make sense.
With this in mind,
we define $D(W)$ to be the subspace of $W^{*}$ consisting of $\alpha$ such that
$$Y^{*}(v,x)\alpha \in W^{*}((x))\;\;\;\;\mbox{ for every }v\in V.$$
Following the proof of Theorem 5.2.1 of [FHL] closely, we easily see 
that FHL actually proved that $D(W)$ is a weak $V$-module.
It is clear that $D(W)$ is the maximal weak $V$-module inside $W^{*}$ with 
$Y=Y^{*}$.
Next, we prove that $(W/C_{2}(W))^{*}\subseteq D(W)$ 
(Proposition \ref{pc2d}). Therefore, any regular vertex operator algebra
satisfies Zhu's $C_{2}$-finiteness 
condition. Furthermore, by exploiting a result of Frenkel and Zhu [FZ]
on fusion rules in terms of $A(V)$-bimodules 
we prove that all fusion rules are finite.
We also define a graded subspace $C_{1}(V)$ of $V$ and prove that any 
graded subspace of $V$ complementary to $C_{1}(V)$ generates $V$ as a vertex 
operator algebra (Proposition \ref{pc1g}). As a corollary we prove that
any regular vertex operator algebra is finitely generated.

In [L5], the notion of $D(W)$ was extended further to the 
notion of ${\cal{D}}(W)$ consisting of what we called
representative functionals on $W$. The space ${\cal{D}}(W)$ was proved to
have a natural $V$-bimodule structure where the $V$-bimodule actions 
were not exactly $Y^{*}$, but certain analytic continuations of $Y^{*}$.
Using the notion ${\cal{D}}(W)$ we were able to prove some results of 
Peter-Weyl type for vertex operator algebras.

This paper is organized as follows: In Section 2, we associate a canonical
weak $V$-module $D(W)$ inside $W^{*}$ to any weak $V$-module $W$
and we prove that $D(W)=W'$ for a certain class $\cal{A}$ of vertex 
operator algebras. In Section 3, we prove the $C_{2}$-finiteness condition
and the finiteness of fusion rules for vertex operator algebras of 
class $\cal{A}$. We also prove that vertex operator algebras of class 
${\cal{A}}$ are finitely generated.

{\bf Acknowledgments.} We would like to thank James Lepowsky 
for useful discussions and suggestions. We would also like to 
thank IAS for its financial support and hospitality
while this work was done during a membership in the Spring of 1997.

\section{The contragredient module and an extension}
In this section, we first review some basic notions and the 
contragredient module [FHL] and we then give a natural extension
of the contragredient module (to a weak module).

We shall use standard definitions and notions in [FHL] and [FLM]
such as the notions of vertex operator algebra, module, intertwining 
operator and fusion rule, which will not be given here.
We shall also use certain concepts which we recall next.

Let $V$ be a vertex operator algebra. A {\em weak} $V$-module [DLM2] is
a pair $(W,Y_{W})$ satisfying all the axioms except
those involving the grading  for a $V$-module given in [FHL] and [FLM]. 
It has been noticed in [DLM2] that the $L(-1)$-derivative property:
$Y_{W}(L(-1)v,x)={d\over dx}Y_{W}(v,x)$ for $v\in V$,  
follows from the other axioms, so that one 
need not check this axiom for having a module or weak module.
A {\em generalized} $V$-module [HL] is a weak $V$ module on which $L(0)$ 
acts semisimply.
An {\em ${\N}$-gradable} weak $V$-module is a weak $V$-module 
$W$ on which there exists 
an ${\N}$-grading $W=\coprod_{n\in {\N}}W(n)$ such that
\begin{eqnarray}
v_{m}W(n)\subseteq W(\wt v+n-m-1)
\end{eqnarray}
for homogeneous $v\in V$ and for $m\in {\Z}, n\in {\N}$, 
where by convention $W(n)=0$ for $n<0$. (This notion was essentially 
introduced by Zhu in [Z].)
It is clear that the sum of ${\N}$-gradable weak $V$-modules is
again an ${\N}$-gradable weak $V$-module.

Let $M=\coprod_{h\in {\C}}M_{(h)}$ be a $V$-module. By definition 
([FHL], [FLM]),
for any $h\in {\C}$, $M_{(h+n)}=0$ for sufficiently small integer $n$.
For any $h\in {\C}$, it is clear that
$\coprod_{m\in {\Z}}M_{(m+h)}$ is a submodule of $M$. Furthermore,
let $k\in {\Z}$ be such that $M_{(h+k)}\ne 0$ and $M_{(h+m)}=0$ for $m<k$.
Then $\coprod_{m\in {\Z}}M_{(m+h)}=\sum_{n\in {\N}}M_{(h+k+n)}$ is an 
${\N}$-gradable $V$-module with $\left(\coprod_{m\in {\Z}}M_{(m+h)}\right)(n)=
M_{(n+k+h)}$ for $n\in {\N}$. It follows that $M$ is a direct sum of
${\N}$-gradable $V$-modules, so that $M$ is an ${\N}$-gradable $V$-module.

A vertex operator algebra $V$ is said to be {\em rational} if
any ${\N}$-gradable weak $V$-module is a direct sum of irreducible 
${\N}$-gradable weak $V$-modules, and $V$ is said to be {\em regular} [DLM2] 
if any weak $V$-module is a direct sum of irreducible (ordinary) 
$V$-modules. If $V$ is rational, it was proved in  [DLM1] that
$V$ has only finitely many inequivalent irreducible modules and that each 
irreducible ${\N}$-gradable weak $V$-module is a module, so that
this notion of rationality is the same as the one defined in [Z].
Note that there are different notions of rationality (see for example [HL]).

Examples of rational vertex operator algebras are: $V_{L}$ 
associated to a positive definite even lattice $L$ ([B], [D1], [FLM]); 
$L(\ell,0)$ associated to 
a finite-dimensional simple Lie algebra $\frak{g}$ and a positive integer 
$\ell$ ([DL], [FZ], [L2]); $L(c_{p,q},0)$ associated to the Virasoro algebra 
and a rational number
$c_{p,q}$ ([FZ], [DMZ], [W]); $V^{\natural}$, Frenkel, Lepowsky and Meurman's 
Moonshine module ([B], [D2], [FLM]); tensor products of 
vertex operator algebras from above. It was proved in [DLM2] that all 
these rational vertex operator algebras are also regular.

It is well known (cf. [B], [FFR], [L1], [L3], [MP]) that there is a 
Lie algebra $g(V)$ associated to any vertex operator 
algebra $V$. More precisely,
\begin{eqnarray}
g(V)=\hat{V}/d\hat{V}
\end{eqnarray}
where 
$$\hat{V}=V\otimes {\C}[t,t^{-1}]\;\;\mbox{and } 
d=L(-1)\otimes 1+1\otimes {d\over dt},$$
with the following bracket formula:
\begin{eqnarray}\label{ecomm}
[u(m),v(n)]=\sum_{i=0}^{\infty}{m\choose i}(u_{i}v)(m+n-i)
\end{eqnarray}
for $u,v\in V, m,n\in {\Z}$, where $u(m)=u\otimes t^{m}+d\hat{V}$.
Furthermore, $g(V)$ is a ${\Z}$-graded Lie algebra where 
$\deg u(m)=\wt u-m-1$ for homogeneous $u$ and for $m\in {\Z}$.
Let $g(V)_{\pm}$ be the subalgebras of $g(V)$ linearly spanned
by homogeneous elements of positive degrees (negative degrees).
A $g(V)$-module $W$ is said to be {\em restricted} if for any $v\in V, w\in W$,
$v(m)w=0$ for $m$ sufficiently large.
It is easy to see that any weak $V$-module $W$ is a restricted 
$g(V)$-module where $v(n)$ for $v\in V, n\in {\Z}$ is represented
by $v_{n}$. 

Let $V$ be a vertex operator algebra, let $W=\coprod_{h\in {\C}}W_{(h)}$ 
be a $V$-module, and
let $W'=\coprod_{h\in {\C}}W_{(h)}^{*}$ be the restricted dual of $W$.
For $v\in V, w'\in W'$, we define
\begin{eqnarray}
\<Y'(v,x)w',w\>
=\left\<w',Y\left(e^{xL(1)}\left(-x^{-2}\right)^{L(0)}v,x^{-1}\right)w
\right\>.
\end{eqnarray}

The following fundamental result was due to Frenkel, 
Huang and Lepowsky ([FHL], Theorem 5.2.1 and Proposition 5.3.1).

\bp{pfhl0}
The pair $(W',Y')$ carries the structure of a $V$-module and $(W'',Y'')=(W,Y)$.
\ep

This module is called the {\em contragredient} module of $W$.
Proposition \ref{pfhl0} is analogous to the fact in 
the classical Lie theory that
for any Lie algebra $g$ and any $g$-module $U$, 
$U^{*}$ is a $g$-module where 
$$(af)(u)=-f(au)\;\;\;\mbox{ for }a\in g, u\in U, f\in U^{*}.$$
In the following we consider a natural extension of
the contragredient module (in general to a weak module). 

Let us start with a weak $V$-module $W$ (without grading). 
For $v\in V, \alpha\in W^{*}$ we define
\begin{eqnarray}\label{e*1}
\<Y^{*}(v,x)\alpha,w\>
=\left\<\alpha,Y\left(e^{xL(1)}\left(-x^{-2}\right)^{L(0)}v,x^{-1}\right)w
\right\>.
\end{eqnarray}
If $W$ is an (ordinary) $V$-module, then $Y^{*}$ extends $Y'$.
It follows from Proposition 5.3.1 of [FHL] that
\begin{eqnarray}\label{e*2}
\<\alpha,Y(v,x)w\>
=\left\<Y^{*}\left(e^{xL(1)}\left(-x^{-2}\right)^{L(0)}v,x^{-1}\right)\alpha,w
\right\>.
\end{eqnarray}

\br{rstar1} {\em Since $e^{xL(1)}(-x^{-2})^{L(0)}v$ is a finite sum and
$$Y(u,x^{-1})w\in W((x^{-1}))\;\;\;\mbox{ for any }u\in V, w\in W,$$
by (\ref{e*1}) we have
\begin{eqnarray}
\<Y^{*}(v,x)\alpha,w\>\in {\C}((x^{-1}))).
\end{eqnarray}
That is,}
\begin{eqnarray}
Y^{*}(v,x)\alpha\in {\rm Hom}(W,{\C}((x^{-1}))).
\end{eqnarray}
\er

For $v\in V$, we set
\begin{eqnarray}
Y^{*}(v,x)=\sum_{n\in {\Z}}v^{*}_{n}x^{-n-1}.
\end{eqnarray}
Let $v\in V_{(h)}$, {\em i.e.}, $L(0)v=hv$. Then
\begin{eqnarray}\label{edual*}
\< v^{*}_{n}\alpha,w\>=\left\<\alpha,
(-1)^{h}\sum_{i\in {\N}}{1\over i!}(L(1)^{i}v)_{2h-n-i-2}w\right\>.
\end{eqnarray}
Furthermore, if $v$ is quasi-primary, {\em i.e.}, $L(1)v=0$, then
\begin{eqnarray}\label{*quasi}
\< v^{*}_{n}\alpha,w\>=\<\alpha,(-1)^{h}v_{2h-n-2}w\>.
\end{eqnarray}

It was observed in [HL] that FHL ([FHL], Proposition \ref{pfhl0})
in fact proves the following opposite Jacobi identity:
\begin{eqnarray}\label{ehl}
& &x_{0}^{-1}\delta\left(\frac{x_{1}-x_{2}}{x_{0}}\right)
Y_{W}(e^{x_{2}L(1)}(-x_{2}^{-2})^{L(0)}v,x_{2}^{-1})
Y_{W}(e^{x_{1}L(1)}(-x_{1}^{-2})^{L(0)}u,x_{1}^{-1})\nonumber\\
& &-x_{0}^{-1}\delta\left(\frac{x_{2}-x_{1}}{-x_{0}}\right)
Y_{W}(e^{x_{1}L(1)}(-x_{1}^{-2})^{L(0)}u,x_{1}^{-1})
Y_{W}(e^{x_{2}L(1)}(-x_{2}^{-2})^{L(0)}v,x_{2}^{-1})\nonumber\\
&=&x_{2}^{-1}\delta\left(\frac{x_{1}-x_{0}}{x_{2}}\right)
Y_{W}(e^{x_{2}L(1)}(-x_{2}^{-2})^{L(0)}Y(u,x_{0})v,x_{2}^{-1})
\end{eqnarray}
for $u,v\in V$. Although this observation was made in [HL] for
a module $W$, obviously this is true if $W$ is a weak module.
(Notice that the symbol $Y^{*}(v,x)$ was used 
in [HL] for $Y(e^{xL(1)}(-x^{-2})^{L(0)}v,x^{-1})$, which 
is {\em different} from ours.)

Taking ${\rm Res}_{x_{0}}$ from (\ref{ehl}) and then using
(\ref{e*1}) we obtain
\begin{eqnarray}
[Y^{*}(u,x_{1}),Y^{*}(v,x_{2})]
={\rm Res}_{x_{0}}x_{2}^{-1}\delta\left(\frac{x_{1}-x_{0}}{x_{2}}\right)
Y^{*}(Y(u,x_{0})v,x_{2}).
\end{eqnarray}
The same proof of Theorem 5.2.1 of [FHL] shows that 
$Y^{*}(L(-1)v,x)=d/dxY^{*}(v,x)$.
Then we have proved:

\bp{plie}
%(a) There is an anti-automorphism $\sigma$ of $g(V)$ such that
%\begin{eqnarray}
%\sigma(v(n))=\sum_{i\in {\N}}{(-1)^{h}\over i!}(L(1)^{i}v)(2h-n-i-2)
%\end{eqnarray}
%for $v\in V_{(h)},h, n\in {\Z}$.
Let $W$ be any weak $V$-module. Then $W^{*}$ is a $g(V)$-module.$\;\;\;\;\Box$
\ep

To obtain a weak $V$-module out of $(W^{*},Y^{*})$ we consider
the Jacobi identity for $Y^{*}$. Notice that in the definition of
a (weak) module, the truncation condition, in this case which is
 $Y^{*}(v,x)\alpha\in W^{*}((x))$, is necessary
for the Jacobi identity to make sense. 
For example, the first term
$$x_{0}^{-1}\delta\left(\frac{x_{1}-x_{2}}{x_{0}}\right)Y^{*}(u,x_{1})
Y^{*}(v,x_{2})\alpha$$
of the Jacobi identity may not exist {\em algebraically} 
if $Y^{*}(v,x_{2})\alpha$
involves infinitely many negative powers of $x_{2}$.
Having known this fact, we consider a certain subspace of $W^{*}$.

\bd{dD}
{\em Let $W$ be a weak $V$-module. Then we define $D(W)$ to be the 
subspace of $W^{*}$ consisting of vectors $\alpha$ 
such that for every $v\in V$
\begin{eqnarray}\label{etrun}
Y^{*}(v,x)\alpha\in W^{*}((x)),
\end{eqnarray}
{\em i.e.}, $\;v^{*}_{n}\alpha=0$ for $n$ sufficiently large.}
\ed

If $W$ is an (ordinary) $V$-module, it is clear that $W'\subseteq D(W)$.
The following proposition gives a characterization of elements of $D(W)$.

\bp{pcdw} Let $W$ be a weak $V$-module and let $\alpha\in W^{*}$.
Then $\alpha\in D(W)$ if and only if for any $v\in V$, there exists 
$k\in {\Z}$  such that
\begin{eqnarray}
v_{m}W\subseteq \ker \alpha \;\;\;\;\mbox{ for }m\le k,
\end{eqnarray}
or equivalently, there exists $r\in {\Z}$  such that
\begin{eqnarray}
x^{r}\<\alpha, Y(v,x)w\>\in {\C}[x^{-1}]
\end{eqnarray}
for all $w\in W$.
\ep

\pf Notice that $e^{xL(1)}(-x^{-2})^{L(0)}v$ for $v\in V$ is a finite sum.
Then it follows from (\ref{e*1}) and (\ref{e*2}) immediately. $\;\;\;\;\Box$

As a corollary we have:

\bc{ccdw} Let $W$ be a weak $V$-module and let $U$ be a subspace of $W$.
Then $(W/U)^{*}\subseteq D(W)$ if and only for any $v\in V$ there exists
$k\in {\Z}$ such that $v_{m}W\subseteq U$ for $m\le k$, where
$(W/U)^{*}$ is viewed as a subspace of $W^{*}$ in the natural way.
$\;\;\;\;\Box$
\ec

The following lemma establishes the stability of the action of $g(V)$ on
$D(W)$.

\bl{lsub} Let $W$ be a weak $V$-module. Then $D(W)$ is a restricted 
$g(V)$-submodule of $W^{*}$.
\el

\pf Since $W^{*}$ is a $g(V)$-module, 
for $u,v\in V, m\in {\Z},\alpha\in D(W)$, we have
\begin{eqnarray}\label{esta}
Y^{*}(u,x)v^{*}_{m}\alpha=v^{*}_{m}Y^{*}(u,x)\alpha
-\sum_{i\in {\N}}{m\choose i}x^{m-i}Y^{*}(v_{i}u,x)\alpha.
\end{eqnarray}
Because $Y^{*}(v_{i}u,x)\alpha\in W^{*}((x))$ for each $i\in {\N}$ and
$v_{i}u=0$ for all but finitely many $i\in {\N}$, we have
$$Y^{*}(u,x)v_{m}^{*}\alpha\in W^{*}((x)).$$
Thus $v^{*}_{m}\alpha\in D(W).$
Therefore $D(W)$ is a $g(V)$-submodule of $W^{*}$. From the definitions, 
$D(W)$ is restricted. $\;\;\;\;\Box$

Furthermore we have:

\bp{pfhl2}
Let $W$ be a weak $V$-module. Then $D(W)$ is a weak $V$-module. 
\ep

\pf Notice that the existences of the three main terms in the Jacobi identity
for $Y^{*}$ are guaranteed by the truncation condition (\ref{etrun}). 
Then it follows from Lemma \ref{lsub} and (\ref{ehl})
immediately.$\;\;\;\;\Box$

\br{rdm} 
{\em  Because of the truncation axiom in the notion of (weak) module, it 
is clear that $D(W)$ is the maximal weak $V$-module in $W^{*}$ with 
$Y=Y^{*}$.}
\er

%\br{rlt} 
%{\em Notice that in the definition of $D(W)$, we require a global truncation:
%$Y^{*}(u,x)\alpha\in W^{*}((x))$ for $u\in V$. That is, for any $w\in W$,
%$\<Y^{*}(u,x)\alpha,w\>$ is a Laurent polynomial where for fixed $u$ 
%and $\alpha$, the orders of the possible pole at $x=0$ are bounded. 
%One may consider $\alpha\in W^{*}$ with local truncation property:
%$$\<Y^{*}(u,x)\alpha,w\>\in {\C}[x,x^{-1}]\;\;\;\;\mbox{ for }
%u\in V, w\in W.$$
%Let $LP(W)$ be the space of all linear functionals with the above 
%local truncation property. Then using (\ref{esta}) one can show that
%$LP(W)$ is stable under the action of $g(V)$.
%But $LP(W)$ is not a (weak) $V$-module because as mentioned before one needs
%the global truncation property for the Jacobi identity to make sense.}
%\er

Let $W$ be a weak $V$-module. For any $h\in {\C}$, we set
\begin{eqnarray}
W_{(h)}=\{w\in W\;|\; L(0)w=hw\}.
\end{eqnarray}
Furthermore we set
$$W^{0}=\coprod_{h\in {\C}}W_{(h)}.$$
Then it is clear that $W^{0}$ is a submodule of $W$.

\bp{pdsub}
Let $W$ be an (ordinary) $V$-module. Then $D(W)^{0}=W'$.
\ep

\pf Let $\alpha\in D(W)_{(h_{1})}, w\in W_{(h_{2})}$. Then
$$h_{1}\<\alpha,w\>=\<L(0)\alpha,w\>=\<\alpha,L(0)w\>=h_{2}\<\alpha,w\>.$$
It follows that 
$$\<D(W)_{(h_{1})},W_{(h_{2})}\>=0\;\;\;\mbox{ for }h_{1}\ne h_{2}.$$
Thus $D(W)_{(h)}\subseteq W_{(h)}^{*}$ for $h\in {\C}$.
Because $W_{(h)}^{*}\subseteq D(W)_{(h)}$ we get
$D(W)_{(h)}=W_{(h)}^{*}$ for $h\in {\C}$. Thus $D(W)^{0}=W'$. $\;\;\;\;\Box$

By definition $D(W)=D(W)^{0}$ if and only if
$L(0)$ acts semisimply on $D(W)$. 
Then as an immediate corollary we have:

\bc{cL(0)}
Let $V$ be a vertex operator algebra such that $L(0)$ acts semisimply 
on any weak module. Then for any $V$-module $W$ we have $D(W)=W'$.
$\;\;\;\;\Box$
\ec

Suppose that $V$ is regular. Then any weak module is a direct sum of
irreducible $V$-modules, so that $L(0)$ acts semisimply.
Then we immediately have:

\bc{creg1}
Let $V$ be a regular vertex operator algebra and let $W$ be a $V$-module. 
Then $D(W)=W'$. $\;\;\;\;\Box$
\ec

Motivated Corollary \ref{cL(0)} we define $\cal{A}$ to be the class of 
vertex operator algebras satisfying the condition that
$L(0)$ acts semisimply on any weak module. 
Let $V$ be a vertex operator algebra containing a vertex operator 
subalgebra (with the same Virasoro element) $V^{0}$ of class $\cal{A}$. 
Since any weak $V$-module is a weak $V^{0}$-module, $L(0)$ (the same 
for both $V$ and $V^{0}$) acts semisimply on any weak $V$-module. Thus
$V$ is also of class $\cal{A}$. Since
the class $\cal{A}$ contains all regular 
vertex operator algebras, any vertex operator algebra that has a regular 
vertex operator subalgebra (with the same 
Virasoro element) is of class $\cal{A}$. In the next section, we shall study
some finiteness properties for this class of vertex operator algebras.

Notice that in the definition of $D(W)$, it was required that 
$Y^{*}(v,x)\alpha\in W^{*}((x))$ for all $v\in V$. However, 
for vertex operator algebras of certain types such as those associated 
to affine Lie algebras or the Virasoro algebra we only need to check 
this for each $v$ of a (usually finite-dimensional) subspace.

\bp{pmin}
Let $V$ be a vertex operator algebra such that $V=\coprod_{n\ge 0}V_{(n)}$
(without negative weights) with $V_{(0)}={\C}{\bf 1}$ and let
$U$ be a graded subspace of $\coprod_{n\ge 1}V_{(n)}$ ($\subseteq V$)
such that $V$ is linearly spanned by elements
\begin{eqnarray}\label{espan}
{\bf 1},\;\;\;u^{1}_{-n_{1}}\cdots u^{r}_{-n_{r}}{\bf 1},
\end{eqnarray}
where $u^{i}\in U, n_{1},\dots, n_{r}\ge 1$.
Let $W$ be a weak $V$-module and let $\alpha\in W^{*}$. If for any $u\in U$,
there exists $k\in {\Z}$ such that for all $w\in W$,
\begin{eqnarray}\label{ered1}
\<\alpha, u_{m}w\>=0\;\;\;\mbox{ whenever }m\le k.
\end{eqnarray}
Then $\alpha\in D(W)$.
\ep

\pf Let $B$ be the subspace of $V$ consisting of each $v$ such that there 
exists $k\in {\Z}$ such that
$$x^{k}\<\alpha,Y(v,x)w\>\in {\C}[x^{-1}]
\;\;\;\mbox{ for all }w\in W.$$
Then by Proposition \ref{pcdw}
$\alpha\in D(W)$ if and only if $V\subseteq B$.
In the following we shall prove by induction that
$$\oplus_{i=0}^{n}V_{(i)}\subseteq B\;\;\;\mbox{ for }n\in {\N}.$$
First, from the assumption we have $U\subseteq B$.
Since $V_{(0)}={\bf C}{\bf 1}$, it is clear that $V_{(0)}\subseteq B$.
Assume that $\oplus_{i=0}^{n}V_{(i)}\subseteq B$ for some $n\in {\N}$.
Let $a=u_{-m}v\in V_{(n+1)}$ be a homogeneous element where 
$u\in U, v\in V, m\ge 1$. 
Note that for any $w\in W$, by taking 
${\rm Res}_{x_{1}}{\rm Res}_{x_{0}}x_{0}^{-m}$ from the Jacobi identity 
we obtain
\begin{eqnarray}
\<\alpha,Y(a,x_{2})w\>
&=&{\rm Res}_{x_{1}}(x_{1}-x_{2})^{-m} 
\<\alpha,Y(u,x_{1})Y(v,x_{2})w\>\nonumber\\
& &-(-x_{2}+x_{1})^{-m}\<\alpha,Y(v,x_{2})Y(u,x_{1})w\>.
\end{eqnarray}
Since $n+1=\wt a=\wt u+m-1+\wt v$ and $\wt u\ge 1$, we have $\wt v\le n$ 
so that 
$v\in B$ (by the inductive assumption). Then there is $k_{1}\in {\Z}$ such that
\begin{eqnarray}\label{ek1}
{\rm Res}_{x_{1}}x_{2}^{k_{1}}(-x_{2}+x_{1})^{-m}
\<\alpha,Y(v,x_{2})Y(u,x_{1})w\>
\in {\C}[x_{2}^{-1}]
\end{eqnarray}
for all $w\in W$. Since $u\in U$, by assumption
there is $r\in {\N}$ such that
$\<\alpha, u_{n}W\>=0$ for $n\le -m-r$.
Then
\begin{eqnarray}\label{esec}
& &{\rm Res}_{x_{1}}(x_{1}-x_{2})^{-m} 
\<\alpha,Y(u,x_{1})Y(v,x_{2})w\>\nonumber\\
&=&\sum_{i=0}^{r}(-1)^{i}{-m\choose i}x_{2}^{i}
\<\alpha,u_{-m-i}Y(v,x_{2})w\>\nonumber\\
&=&\sum_{i=0}^{r}(-1)^{i}{-m\choose i}x_{2}^{i}
\<\alpha,Y(v,x_{2})u_{-m-i}w\>\nonumber\\
& &+\sum_{i=0}^{r}\sum_{j\in {\N}}(-1)^{i}{-m\choose i}{-m-i\choose j}
x_{2}^{-m-j}\<\alpha,Y(u_{j}v,x_{2})w\>.
\end{eqnarray}
Since $\wt u_{j}v=\wt u+\wt v-j-1<\wt u+\wt v+m-1=n+1$, $u_{j}v\in B$ 
for $j\in {\N}$. Because $v,\; u_{j}v\in B$ and $u_{j}v\ne 0$ only 
for finitely many $j\in {\N}$, by (\ref{esec}) there is $k_{2}\in {\Z}$ 
such that
\begin{eqnarray}\label{ek2}
{\rm Res}_{x_{1}}x_{2}^{k_{2}}(x_{1}-x_{2})^{-m} 
\<\alpha,Y(u,x_{1})Y(v,x_{2})w\>\in {\C}[x_{2}^{-1}]
\end{eqnarray}
for all $w\in W$. Combining (\ref{ek1}) with (\ref{ek2}) we get
$$x_{2}^{k}\<\alpha,Y(a,x_{2})w\>\in {\C}[x_{2}^{-1}]$$
for all $w\in W$, where $k=\min \{k_{1},k_{2}\}$. Thus $a\in B$. 
By the spanning property (\ref{espan}), 
$V_{(n+1)}$ is linearly spanned by elements like
$a$. Thus $V_{(n+1)}\subseteq B$. Therefore $V\subseteq B$.
This proves that $\alpha\in D(W)$. $\;\;\;\;\Box$

Proposition \ref{pmin} will be useful if one wants to determine
$D(W)$ explicitly for certain vertex operator algebras.

\section{The $C_{2}$-finiteness condition and the finiteness of fusion rules}
This section is the core of the paper. 
In this section we define subspaces $C_{n}(W)$ for $n\ge 1$ and 
for any weak $V$-module,
generalizing  Zhu's $C_{2}(V)$ subspace defined in [Z]. We
prove that $(W/C_{n}(W))^{*}\subseteq D(W)$ and that 
$V$ is finitely generated if $C_{1}(V)$ is finite-codimensional.
By applying Corollaries \ref{cL(0)} and \ref{creg1}
we then prove the $C_{2}$-finiteness condition, the finiteness of fusion rules
and the finite generating property for the class $\cal{A}$ of vertex 
operator algebras defined in Section 2.

Let $W$ be a weak $V$-module and let $n\ge 2$. We define
$C_{n}(W)$ to be the linear span of elements of type
\begin{eqnarray}\label{ecn}
v_{-n}w,\;\;\;\mbox{ for }v\in V, w\in W.
\end{eqnarray}
Since $(L(-1)v)_{-n}=nv_{-n-1}$, we have:
$$v_{-m}w\in C_{n}(W)\;\;\;\mbox{ for }v\in V, w\in W, m\ge n.$$
Thus
$$\cdots \subseteq C_{n+1}(W)\subseteq C_{n}(W)\subseteq \cdots 
\subseteq C_{2}(W).$$
Note that $C_{2}(V)$ is exactly the one defined in [Z]. 
A vertex operator algebra
$V$ is said to satisfy {\em $C_{2}$-finiteness condition} ([Z], [DLM3]) if 
$V/C_{2}(V)$ is finite-dimensional. (Note that Zhu defined the concept of
finiteness condition $C$ by requiring that $V/C_{2}(V)$ is 
finite-dimensional and that $V$ is a sum of lowest weight modules for
the Virasoro algebra.)

Let $W$ be a generalized $V$-module, {\em i.e.}, a weak module on which
$L(0)$ acts semisimply. Set
\begin{eqnarray}
W_{+}=\oplus_{h\in {\C}, {\rm Re}\;h>0}W_{(h)}.
\end{eqnarray}
Then we define $C_{1}(W)$ to be the 
subspace of $W$ linearly spanned by 
\begin{eqnarray}
v_{-1}W_{+},\;\;\; L(-1)W\;\;\;\mbox{  for  }v\in V_{+}.
\end{eqnarray}
(The subspace $C_{1}(W)$ is designed to make the proof
of Proposition \ref{pc1g} work. However, one may define $C_{1}(W)$ 
differently for other purposes.)
Notice that
if we had defined $C_{1}(V)$ to be the linear span of elements of 
type (\ref{ecn}) with $n=1$, then we would have
$C_{1}(V)=V$ because $v=v_{-1}{\bf 1}\in C_{1}(V)$ for $v\in V$.
That is not what we want.

Since $k!v_{-k-1}=(L(-1)^{k}v)_{-1}$ for $v\in V, k\in {\N}$, we have
\begin{eqnarray}\label{ec1k}
v_{-k-1}w\in C_{1}(W)\;\;\;\mbox{ for }v\in V_{+}, w\in W_{+}, k\in {\N}.
\end{eqnarray}
More generally, we have:

\bl{lc1c} Let $W$ be a generalized $V$-module, let $v\in V, w\in W$ 
be homogeneous and let $r,\;s\in {\N}$ be such that
$L(-1)^{r}v\in V_{+}$ and $L(-1)^{s}w\in W_{+}$.
Then 
\begin{eqnarray}
v_{-r-s-k}w\in C_{1}(W)\;\;\;\mbox{ for }k\ge 1.
\end{eqnarray}
\el

\pf Since $L(-1)^{r}v\in V_{+},\;L(-1)^{s}w\in W_{+}$ by (\ref{ec1k}) we have
$$(L(-1)^{r}v)_{-k}L(-1)^{s}w\in C_{1}(W)\;\;\;\mbox{ for }k\ge 1.$$
Since $[L(-1),u_{m}]=-mu_{m-1}$ for $u\in V, m\in {\Z}$ and
 $L(-1)W\subseteq C_{1}(W)$, we get
$$(L(-1)^{r}v)_{-k-s}w\in C_{1}(W).$$
Then
$$v_{-r-s-k}w\in C_{1}(W)\;\;\;\mbox{ for }k\ge 1.\;\;\;\;\Box$$

\br{rc12} {\em 
In general, $C_{2}(W)$ is not a subspace of $C_{1}(W)$. However, if 
$V=\coprod_{n\in {\N}}V_{(n)}$ (without negative weights)
such that $V_{(0)}={\C}{\bf 1}$, then $C_{2}(V)\subseteq C_{1}(V)$
because $v_{-2}{\bf 1}=L(-1)v$ and ${\bf 1}_{-2}v=0$ for any $v\in V$.}
\er

The following proposition gives a way to find a relatively small
generating subspace of $V$ as a vertex operator algebra by using $C_{1}(V)$.

\bp{pc1g}
Let $V$ be a vertex operator algebra
and $U$ be a graded subspace of $V$ such that $V=U+ C_{1}(V)$ 
and $\oplus_{n\le 0}V_{(n)}\subseteq U+{\C}{\bf 1}$. Then
$V$ is linearly spanned by elements of type
\begin{eqnarray}\label{pbw}
u^{1}_{n_{1}}\cdots u^{r}_{n_{r}}{\bf 1},
\end{eqnarray}
where $r\in {\N}, u^{i}\in U, n_{i}\in {\Z}$ for $1\le i\le r$, {\em i.e.},
$U$ generates $V$ as a vertex operator algebra.
\ep

\pf Let $\<U\>$ be the vertex operator subalgebra of $V$ generated by
$U$. We shall prove by induction that for any $n\in {\N}$,
$$\oplus_{i\le n}V_{(i)}\subseteq \<U\>.$$
By the assumption, we have 
$\oplus_{i\le 0}V_{(i)}\subseteq U+{\C}{\bf 1}\subseteq \<U\>$.
Suppose that $\oplus_{i\le n}V_{(i)}\subseteq \<U\>$ for some $n\in {\N}$.
Let $a\in V_{(n+1)}$. Since 
$V_{(n+1)}=C_{1}(V)\cap V_{(n+1)}+U\cap V_{(n+1)}$, we have
\begin{eqnarray}
a=u^{1}_{-1}v^{1}+\cdots +u^{r}_{-1}v^{r}+L(-1)u+b
\end{eqnarray}
for some homogeneous $u^{i}, v^{i}, u\in V, b\in U$ with 
$\wt u^{i}_{-1}v^{i}=\wt u^{i}+\wt v^{i}=n+1$, $\wt u^{i}>0$, $\wt v^{i}>0$ 
and $\wt u=n$.
Then $\wt u^{i},\wt v^{i}\le n$ for all $i$. By inductive assumption, we have
$$u^{i},\;\;v^{i},\;\; u\in \<U\>\;\;\;\mbox{ for all }i,$$
so that
$$u^{i}_{-1}v^{i}\in \<U\>\;\;\;\mbox{ for all }i$$ 
and $L(-1)u=u_{-2}{\bf 1}\in \<U\>$. Thus 
$a\in \<U\>$. This proves that $V_{(n+1)}\subseteq \<U\>$. Therefore 
$V=\<U\>$.$\;\;\;\;\Box$

As a refinement of Proposition \ref{pc1g}, it was proved in [KL] 
that $V$ is linearly spanned by elements in (\ref{pbw}) with
$n_{1},\dots, n_{r}<0$ in a fixed lexicographical order.

As an immediate corollary of Proposition \ref{pc1g} we have:

\bc{cfinite}
Let $V$ be a vertex operator algebra such that $C_{1}(V)$ is
finite-codimensional. Then $V$ is finitely generated.$\;\;\;\;\Box$
\ec

Furthermore, by Remark \ref{rc12} we have:

\bc{cc2finite}
Suppose that $V=\coprod_{n\in {\N}}V_{(n)}$ (without negative weights) 
and $V_{(0)}={\C}{\bf 1}$. Then the $C_{2}$-finiteness
condition on $V$ implies that $V$ is finitely generated.
$\;\;\;\;\Box$
\ec

In the following we shall prove that the $C_{2}$-finiteness condition holds
for vertex operator algebras of class $\cal{A}$. For convenience, we 
set $Z_{n}(W)=W/C_{n}(W)$. 

\bp{pc2d}
Let $W$ be a weak $V$-module and let $n\ge 2$. Then
\begin{eqnarray}
Z_{n}(W)^{*}=\{\alpha\in W^{*}\;|\; v^{*}_{m}\alpha=0\;\;\;\mbox{ for 
homogeneous }v\in V,
m\ge 2\wt v+n-2\}.
\end{eqnarray}
In particular,
$$Z_{n}(W)^{*}\subseteq D(W).$$
\ep

\pf Let $\alpha\in Z_{n}(W)^{*}, v\in V$. Then
\begin{eqnarray}
\<\alpha, v_{-m}w\>=0\;\;\;\;\mbox{ for }v\in V, w\in W, m\ge n.
\end{eqnarray}
Let $u\in V$ be homogeneous. Then for any $m\ge 2\wt u+n-2$
\begin{eqnarray}
\<u^{*}_{m}\alpha,w\>=\left\<\alpha,
\sum_{i\in {\N}}{(-1)^{\wt u}\over i!}(L(1)^{i}u)_{2\wt u-i-m-2}w\right\>=0
\end{eqnarray}
for any $w\in W$. Thus
\begin{eqnarray}\label{eczn}
u^{*}_{m}\alpha=0 \;\;\;\;\mbox{ for }m\ge 2\wt u+n-2.
\end{eqnarray}
Conversely, let $\alpha\in W^{*}$ be such that (\ref{eczn}) holds and 
let $v\in V_{(k)}$ for $k\in {\Z}$. Then
\begin{eqnarray}
\<\alpha, v_{-m}w\>
=\sum_{i\in {\N}}(-1)^{k}{1\over i!}
\left\<(L(1)^{i}v)^{*}_{2k-i-2+m}\alpha,w\right\>=0
\end{eqnarray}
for $m\ge n$ (noticing that $2k-i-2+m=2\wt (L(1)^{i}v)+i+m$). 
Thus $\alpha\in Z_{n}(W)^{*}$. Then the proof is complete.
$\;\;\;\;\Box$

%For a weak $V$-module $W$, let $Z^{*}(W)=\cup_{n\in {\N}}Z_{n}(W)^{*}$.
%\bl{lz*}
%Let $W$ be a weak $V$-module. Then $Z^{*}(W)$ is a submodule of $D(W)$.
%\el

Now we present our key result of this paper:

\bp{pc2}
Let $V$ be a vertex operator algebra of class $\cal{A}$ and let $W$ be a
$V$-module. Then
$Z_{n}(W)$ is finite-dimensional for $n\ge 2$ and $Z_{1}(W)$ is
also finite-dimensional if the real parts of the weights of $W$ are bounded
{}from below. 
\ep

\pf We first notice that a graded subspace $U$ of $W$ is 
finite-codimensional if and only if $(W/U)^{*}\subseteq W'$, where $(W/U)^{*}$
is viewed as a natural subspace of $W^{*}$. 
It is clear that $C_{n}(W)$ is a graded subspace. Then
$\dim Z_{n}(W)<\infty$ if and only if $Z_{n}(W)^{*}\subseteq W'$.
Furthermore, since $D(W)=W'$ (Corollary \ref{cL(0)}),
$\dim Z_{n}(W)<\infty$ if and only if 
$$Z_{n}(W)^{*}\subseteq D(W).$$
If $n\ge 2$, it follows immediately from Proposition \ref{pc2d}.
For $n=1$, since the real parts of the weights of $W$ are bounded from 
below, there 
is a nonnegative integer $s$ such that $L(-1)^{s}W\subseteq W_{+}$.
Then it follows from Lemma \ref{lc1c} and Corollary \ref{ccdw} that
$Z_{1}(W)^{*}\subseteq D(W)$. This completes the proof.
$\;\;\;\;\Box$

Combining Corollary \ref{cfinite} with
Proposition \ref{pc2}, we immediately have:

\bt{tfinite}
Any vertex operator algebra of class $\cal{A}$ satisfies 
the $C_{2}$-finiteness condition and it is finitely generated. 
In particular, any regular vertex operator algebra $V$ 
satisfies the $C_{2}$-finiteness condition and it
is finitely generated. $\;\;\;\;\Box$
\et

Next we shall prove the finiteness of fusion rules for vertex operator algebras
of class $\cal{A}$.
To do this we shall use Frenkel and Zhu's $A(V)$-bimodule theory [FZ].

Recall $A(V)$ and $A(W)$ from [FZ] and [Z].  
For any weak $V$-module $W$, let $O(W)$ be the linear span of elements of type
\begin{eqnarray}
\sum_{i\in {\N}}{\wt v\choose i}v_{i-2}w\;\;\left(={\rm Res}_{x}
\frac{(1+x)^{\wt v}}{x^{2}}
Y(v,x)w\right)
\end{eqnarray}
for homogeneous $v\in V$ and for $w\in W$. It was proved ([Z], Lemma 2.1.1) 
that
\begin{eqnarray}
\sum_{i\in {\N}}{\wt v\choose i}v_{i-m}w\in O(W)
\end{eqnarray}
for homogeneous $v\in V$ and for $w\in W, m\ge 2$.
Set $A(W)=W/O(W)$. Then we have [Z]:

\bp{pZhu}
(a) The space $A(V)$ is an associative algebra with the product:
\begin{eqnarray}
(u+O(V))(v+O(V))=\sum_{i\in {\N}}{\wt u\choose i}(u_{i-1}v+O(V))
\end{eqnarray}
for homogeneous $u,v\in V$.

(b) For any ${\N}$-gradable weak $V$-module $W=\coprod_{n\in {\N}}W(n)$, 
$W(0)$ is an $A(V)$-module
where $v+O(V)$ acts on $W(0)$ as $v_{\wt v-1}$ for homogeneous $v\in V$.

(c) There is a one-to-one correspondence between the set of equivalence 
classes of 
irreducible $A(V)$-modules and the set of equivalence classes of irreducible
${\N}$-gradable weak $V$-modules.
\ep

Furthermore, we have [FZ]:

\bp{pFZ}
Let $W$ be a weak $V$-module. Then $A(W)$ is an $A(V)$-bimodule with
the following left and right actions:
\begin{eqnarray}
& &(v+O(V))(w+O(W))=\sum_{i\in {\N}}{\wt v\choose i}(v_{i-1}w+O(W))
\label{eleft}\\
& &(w+O(W))(v+O(V))=\sum_{i\in {\N}}{\wt v-1\choose i}(v_{i-1}w+O(W))
\label{eright}
\end{eqnarray}
for homogeneous $v\in V$ and for $w\in W$.
\ep
 
Let $W_{1},W_{2},W_{3}$ be irreducible $V$-modules and let 
$I{W_{3}\choose W_{1}W_{2}}$
be the space of intertwining operators of the indicated type. Then we have
([L4], Proposition 2.10, and [FZ]):

\bp{pfz}
Let $W_{1},W_{2},W_{3}$ be irreducible $V$-modules. Then
\begin{eqnarray}
\dim I{W_{3}\choose W_{1}W_{2}}\le 
\dim {\rm Hom}_{A(V)}(A(W_{1})\otimes _{A(V)}W_{2}(0), W_{3}(0)),
\end{eqnarray}
where $W_{2}(0)$ and $W_{3}(0)$ are the lowest weight subspaces of 
$W_{2}$ and $W_{3}$, respectively.
\ep

The following lemma was proved in 
[DLM3] (Proposition 3.6) (see also [Z], Lemma 4.4.1).

\bl{lz}
Let $W$ be a $V$-module such that $\dim Z_{2}(W)<\infty$. Then
$\dim A(W)<\infty$. 
\el

Combining Lemma \ref{lz} with Propositions \ref{pZhu} and \ref{pc2} we obtain

\bt{tc2f}
Let $V$ be a vertex operator algebra of class $\cal{A}$.
Then there are only finitely many inequivalent irreducible 
${\N}$-gradable weak $V$-modules and 
the fusion rule for any three 
irreducible modules is finite. In particular, the assertions hold 
if $V$ is regular.$\;\;\;\;\Box$
\et

\br{ravd}
{\em Motivated by Proposition \ref{pc2d},
one may also consider the subspace $A(W)^{*}$ of $W^{*}$ for a (weak) 
$V$-module $W$. By using the proof of Proposition \ref{pc2d} for $A(W)^{*}$,
one can see that in general $A(W)^{*}$ may not be a subspace of $D(W)$. 
However, it was proved in [L5] that
$A(W)^{*}$ is a subspace of a canonical weak $V$-bimodule
${\cal{D}}(W)$ where ${\cal{D}}(W)$ is the space of what we call
representative functionals on $W$ containing $D(W)$ as a subspace.}
\er

\br{ram}
{\em Note that $\dim A(W_{1})<\infty$ is a sufficient condition for
the fusion rule $\dim I{W_{3}\choose W_{1}W_{2}}$ to be finite.
However, by no means it is necessary. As a matter of fact,
if $A(W_{1})$ is a finitely generated $A(V)$-bimodule,
the fusion rule $\dim I{W_{3}\choose W_{1}W_{2}}$ is finite even though
$A(W_{1})$ may be infinite-dimensional. Indeed, if $S$ is a 
finite-dimensional subspace of $A(W_{1})$ which generates $A(W_{1})$ as an 
$A(V)$-bimodule, then 
$$\dim\; {\rm Hom}_{A(V)}(A(W_{1})\otimes_{A(V)} W_{2}(0),W_{3}(0))
\le \dim\; {\rm Hom}_{\C}(S\otimes W_{2}(0),W_{3}(0)).$$
By Proposition \ref{pfz}, we get $\dim\;I{W_{3}\choose W_{1}W_{2}}<\infty$.}
\er

For the rest of this section we shall give a sufficient condition
for $A(W)$ to be a finitely generated $A(V)$-bimodule.
For convenience, we assume that $V=\coprod_{n\in {\N}}V_{(n)}$ 
with $V_{(0)}={\C}{\bf 1}$. For any $V$-module $W$, we define
$B(W)$ to be the subspace of $W$ linearly spanned by
\begin{eqnarray}
& &u_{-1}W\;\;\;\;\mbox{ for }u\in V_{+},\\
& &v_{0}W\;\;\;\;\mbox{ for homogeneous }v \;\;\mbox{ with }\wt v\ge 2.
\end{eqnarray}
That is, $B(W)=g(V)_{+}W$, where $g(V)_{+}$ is the subalgebra of $g(V)$
linearly spanned by homogeneous elements of positive degrees.

\bp{pfinite}
Let $W$ be a $V$-module such that $W=\coprod_{n\in {\N}}W_{(h+n)}$ 
for some $h\in {\C}$ and let $W^{0}$ be a graded subspace of $W$ 
such that $W=W^{0}+B(W)$. Then $(W^{0}+O(W))/O(W)$ generates $A(W)$ 
as an $A(V)$-bimodule. In particular, if $\dim\;W/B(W)<\infty$,
$A(W)$ is a finitely generated $A(V)$-bimodule.
\ep

\pf Let $E$ be the $A(V)$-bimodule generated by $W^{0}+O(W)$.
We shall prove by induction that 
$W_{(h+n)}+O(W)\subseteq E$ for $n\in {\N}$.
Since $W^{0}$ and $B(W)$ are graded and $B(W)\cap W_{(h)}=0$, we have
$W_{(h)}\subseteq W^{0}$.
Suppose that 
$$(\oplus _{i=0}^{n}W_{(n+h)})+O(W)\subseteq E$$ 
for some $n\in {\N}$.
Let $w\in W_{(h+n+1)}$. Then 
\begin{eqnarray}
w=u^{1}_{-1}w^{1}+\cdots +u^{r}_{-1}w^{r}+v^{1}_{0}w^{r+1}
+\cdots v^{s}_{0}w^{r+s}+w'
\end{eqnarray}
for some homogeneous $u^{i}, v^{j}\in V_{+}, w^{k}\in W$ such that
$\wt u^{i}\ge 1, \;\wt v^{j}\ge 2$. Since 
$\wt u^{i}_{-1}$, $\wt v_{0}^{j}\ge 1$, we have 
$w^{k}\in \oplus_{t=0}^{n}W_{(t+h)}$, so that
by the inductive assumption we have $w^{k}+O(W)\in E$ for 
$1\le k\le r+s$. Then (recall (\ref{eleft}) and (\ref{eright}))
\begin{eqnarray}
& &(u^{i}+O(V))*(w^{i}+O(W))
=\sum_{p=0}^{\wt u^{i}}{\wt u^{i}\choose p}u^{i}_{p-1}w^{i}+O(W)
\in E,\\
& &(v^{j}+O(V))*(w^{j+r}+O(W))-(w^{j+r}+O(W))*(v^{j}+O(V))\nonumber\\
&=&\sum_{p=0}^{\wt v^{i}-1}{\wt v^{j}-1\choose p}
v^{j}_{p}w^{i+r}+O(W)\in E.
\end{eqnarray}
Since $u^{i}_{p-1}w^{i}\in \oplus_{t=0}^{n}W_{(h+t)}$, by the inductive 
assumption
we have $u^{i}_{p-1}w^{i}+O(W)\in E$ for $p\ge 1, 1\le i\le r$.
Then
$$u^{i}_{-1}w^{i}+O(W)=u^{i}*w^{i}-\sum_{p=1}^{\wt u^{i}}{\wt u^{i}\choose p}
u^{i}_{p-1}w^{i}+O(W)\in E.$$
Similarly, we have $v^{j}_{0}w^{r+j}+O(W)\in E$ for $1\le j\le s$.
Then $w+O(W)\in E$. This completes the induction. Therefore $E=A(W)$.
$\;\;\;\;\Box$

As an immediate corollary of Proposition \ref{pfinite}
and Remark \ref{ram}, we have:

\bc{cram}
Let $W_{1}, W_{2}, W_{3}$ be irreducible $V$-modules such that 
$\dim\; W_{1}/B(W_{1})<\infty$. Then 
$\dim\; I{W_{3}\choose W_{1}W_{2}}<\infty$.$\;\;\;\;\Box$
\ec

\end{document}